\documentclass[10pt]{amsart}
\usepackage{amssymb,latexsym,amsthm}

\theoremstyle{plain}

\newtheorem{lemma}{Lemma}

\newtheorem* {theorem A}{Theorem A}
\newtheorem* {theorem B}{Theorem B}
\newtheorem* {theorem C}{Theorem C}
\newtheorem* {theorem D}{Theorem D}
\newtheorem* {theorem C'}{Theorem C'}
\newtheorem* {corollary A}{Corollary A}
\newtheorem* {corollary A'}{Corollary A'}
\newtheorem* {corollary B'}{Corollary B'}
\newtheorem* {proposition A}{Proposition A}
\newtheorem* {proposition A'}{Proposition A'}

\theoremstyle{definition}
\newtheorem{definition}{Definition}

\newtheorem*{question}{Question}

\theoremstyle{remark}

\def\R{\mathbb{R}}
\def\Z{\mathbb{Z}}
\def\N{\mathbb{N}}
\def\NN{\mathbb{N}}
\def\Q{\mathbb{Q}}
\def\xis{\mathfrak{X}^1(M)}
\def\HR{\mathit{HR}}
\def\RR{\mathcal{R}}
\def\eps{\varepsilon}
\def\xis{\mathfrak{X}^1(M)}
\def\xisr{\mathfrak{X}^r(M)}
\def\diff{\mathrm{Diff}^1(S)}

\begin{document}

\title[Topological mixing for flows]
{Robust transitivity and topological mixing for $C^1$-flows}
\author[Abdenur, Avila and Bochi]{Flavio Abdenur, Artur Avila and
Jairo Bochi}
\date{\today}
\thanks{F.A. was supported by FAPERJ and Prodoc/CAPES, A.A. was
supported by FAPERJ and CNPq, and J.B. was supported by
Profix/CNPq.}

\begin{abstract}
We prove that non-trivial homoclinic classes of $C^r$-generic
flows are topologically mixing. This implies that given $\Lambda$
a non-trivial $C^1$-robustly transitive set of a vector field $X$,
there is a $C^1$-perturbation $Y$ of $X$ such that the
continuation $\Lambda_Y$ of $\Lambda$ is a topologically mixing
set for $Y$. In particular, robustly transitive flows become
topologically mixing after $C^1$-perturbations. These results
generalize a theorem by Bowen on the basic sets of generic Axiom A
flows. We also show that the set of flows whose non-trivial
homoclinic classes are topologically mixing is \emph{not} open and
dense, in general.
\end{abstract}

\maketitle



\section{Statement of the Results}

Throughout this paper $M$ denotes a compact $d$-dimensional
boundaryless manifold, $d \geq 3$, and $\xisr$ is the space of
$C^r$ vector fields on $M$ endowed with the usual $C^r$ topology,
where $r \geq 1$ . We shall prove that, generically (residually)
in $\xisr$, nontrivial homoclinic classes are topologically
mixing. As a consequence, nontrivial $C^1$-robustly transitive
sets (and $C^1$-robustly transitive flows in particular) become
topologically mixing after arbitrarily small $C^1$-perturbations
of the flow.

These results generalize the following theorem by Bowen \cite{B}:
non-trivial basic sets of $C^r$-generic Axiom A flows are
topologically mixing. Note that $C^1$-robustly transitive sets are
a natural generalization of hyperbolic basic sets; they are the
subject of several recent papers, such as \cite{BD1} and
\cite{BDP}.

In order to announce precisely our results, let us introduce some
notations and definitions.

Given $t \in \R$ and $X \in \xisr$, we shall denote by $X^t$ the
induced time $t$ map. A subset $\RR$ of $\xisr$ is \emph{residual}
if it contains the intersection of a countable number of open
dense subsets of $\xisr$. Residual subsets of $\xisr$ are dense.
Given an open subset $U$ of $\xisr$, then property (P) is
\emph{generic in $U$} if it holds for all flows in a residual
subset $\RR$ of $U$; (P) is \emph{generic} if it is generic in all
of $\xisr$.

A compact invariant set for $X$ is \emph{non-trivial} if it is
neither a periodic orbit nor a single point. A compact invariant
set $\Lambda$ of $X$ is \emph{transitive} if it there is some
point $x \in \Lambda$ such that the future orbit $\{X^t(x): t
> 0\}$ of $x$ is dense in $\Lambda$; $\Lambda$ is \emph{topologically
mixing} for $X$ if given any nonempty open subsets $U$ and $V$ of
$\Lambda$ then there is some $t_0 > 0$ such that $X^t(U) \cap V
\neq \varnothing$ for all $t \geq t_0$. A non-trivial
$X$-invariant transitive set $\Lambda$ is \emph{$\Omega$-isolated}
if there is some open neighborhood $U$ of $\Lambda$ such that $U
\cap \Omega(X) = \Lambda$. Furthermore, $\Lambda$ is
\emph{isolated} if there is a neighborhood $U$ of $\Lambda$
(called an \emph{isolating block}) such that $$
\Lambda=\bigcap_{t\in\R} X^t(U). $$

Given a hyperbolic closed orbit $\gamma$ of $X$, the
\emph{homoclinic class} of $\gamma$ relative to $X$ is given by $$
H_X(\gamma) = \overline{W^s(\gamma) \pitchfork W^u(\gamma)}, $$
where $\pitchfork$ denotes points of transverse intersection of
the invariant manifolds. $H_X(\gamma)$ is a transitive compact
$X$-invariant subset of the non-wandering set $\Omega(X)$.
Moreover, if $\gamma$ is a closed orbit of index $i$, then the set
$P_i(H_X(\gamma)) \equiv \{p \in H_X(\gamma)\cap \text{Per}(X): p
\text{ is hyperbolic with index $i$} \}$ is dense in
$H_X(\gamma)$.
(See~\cite{BDP}). $H_X(\gamma)$ is not necessarily hyperbolic, but
if $X$ is Axiom A then its basic sets are hyperbolic homoclinic
classes. In the absence of ambiguity, we may write $H(\gamma)$ for
$H_X(\gamma)$.

An \emph{attractor} is a transitive set $\Lambda$ of $X$ that
admits a neighborhood $U$ such that $$ X^t(U) \subset U \text{ for
all } t > 0 \text{, and } \bigcap_{t \in \R} X^t(U) = \Lambda. $$
A \emph{repeller} is an attractor for $-X$. Clearly any attractor
or repeller is $\Omega$-isolated.

An isolated $X$-invariant compact set $\Lambda$ is
\emph{$C^1$-robustly transitive} if there is some open
neighborhood $\mathcal{V}$ of $X$ in $\xis$ and some isolating
block $U$ of $\Lambda$ such that given any $Y \in \mathcal{V}$,
then $$ \Lambda_Y \equiv \bigcap_{t\in\R} Y_t(U) $$ is a compact
transitive non-trivial set of $Y$.

Finally, a flow $X$ is \emph{$C^1$-robustly transitive} if there
is some open neighborhood $\mathcal{W}$ of $X$ in $\xis$ such that
given any $Y \in \mathcal{W}$ then $Y$ is transitive.

\medskip

Our main result is the following:

\begin{theorem A}
There is a residual subset $\RR$ of $\xis$ such that if $\Lambda$
is an isolated non-trivial transitive set of $X \in \RR$, then
$\Lambda$ is topologically mixing for $X$.
\end{theorem A}

Theorem A has the following immediate consequence for robustly
transitive sets or flows:

\begin{corollary A'}
Let $\Lambda$ be a non-trivial robustly transitive set, with
$\mathcal{V}$ and $U$ as in the definition above. Then there is
some residual subset $\RR$ of $\mathcal{V}$ such that if $Y \in
\RR$ then $\Lambda_Y$ is topologically mixing for $Y$. In
particular, given an open set $\mathcal{W} \subset \xis$ of
transitive flows, then there is some residual subset $\RR$ of
$\mathcal{W}$ such that any $Y \in \RR$ is topologically mixing.
\end{corollary A'}

Theorem A is very much a nonhyperbolic, $C^1$ version of Bowen's
aforementioned result. It is actually a consequence of the proof
of the following result:

\begin{theorem B}
Given any $r \in \NN$, there is a residual subset $\RR$ of $\xisr$
such that if $Y \in \RR$ and $H(\gamma)$ is a non-trivial
homoclinic class of $Y$, then $H(\gamma)$ is topologically mixing
for $Y$.
\end{theorem B}

Theorem B follows from general properties of homoclinic classes
combined with simple topological arguments. All of the arguments
in the proof of Theorem B hold in any $C^r$ topology with $r \geq
1$, whereas Theorem A requires the use of $C^1$-generic properties
which are not known in finer topologies.

Pugh's General Density Theorem [Pu] and Theorem B of~\cite{BD2}
(which is stated for diffeomorphisms but holds for flows via the
same arguments) imply that $C^1$-generically any $\Omega$-isolated
transitive set coincides with some homoclinic class. Therefore
Theorem~B implies the following corollary:

\begin{corollary B'}
There is a residual subset $\RR$ of $\xis$ such that if $Y \in
\RR$ and $\Lambda$ is a non-trivial transitive $\Omega$-isolated
set of $Y$, then $\Lambda$ is topologically mixing for $Y$. In
particular, $C^1$-generically any non-trivial attractor/repeller
is topologically mixing.
\end{corollary B'}

Note that Corollary B' generalizes the ``mixing'' aspect of
\cite{MP}. We remark that the dependence of our proofs on the
($C^1$) Closing and Connecting Lemmas means that extending our
results (with the exception of Theorem B) to finer topologies is
probably very difficult.

Of course, in general not every non-trivial homoclinic class is
topologically mixing: the basic sets of suspensions of Axiom A
diffeomorphisms, for example, are not mixing. One may therefore
ask how large is the set of the flows that have a non-trivial
homoclinic class which is \emph{not} mixing. A partial answer to
this question is given by:

\begin{theorem C}
There exists a $4$-manifold $M$ and an open set $\mathcal{U}
\subset \xis$ such that each flow $X$ in a dense subset
$\mathcal{D} \subset \mathcal{U}$ has a non-trivial homoclinic
class which is~\emph{not} topologically mixing for $X$.
\end{theorem C}

Theorem C shows that the residual set $\RR$ of Theorem B is, in
general, not open. The construction in Theorem C relies on the
wild diffeomorphisms from~\cite{BD2} and \cite{BD3}.

On the other hand, robustly transitive flows have relatively tame
dynamics. We pose the following:

\begin{question}
Is the set of ($C^1$-)robustly topologically mixing flows dense in
the set of robustly transitive flows?
\end{question}

In~\cite{AA} the first two authors prove analogues of Theorems A,
B, and C for diffeomorphisms. In addition, a robustly transitive
but non-mixing diffeomorphism is constructed.

The next section first lists some definitions and properties
needed for the proofs and then sets out the proofs themselves.

\section{The Proofs}

Given a hyperbolic periodic point $p$, let $\gamma = \gamma(p)$ be
its orbit and $\Pi_X(p) = \Pi_X(\gamma)$ be its period. Set also
\begin{alignat*}{2}
W^s(p)            &= \{ x \in M:         &&\  d(X^t(x), \gamma(p))
\to 0 \text{ as $t \to + \infty$} \},  \\ W^{ss}(p)         &= \{
x \in M:         &&\  d(X^t(x), X^t(p)) \to 0 \text{ as $t \to +
\infty$} \}.
\end{alignat*}
We define $W^u(p)$ and $W^{uu}(p)$ as the corresponding sets for
$-X$. Note that the set $W^s(p)$ is $X^t$-invariant for all $t \in
\R$, whereas $W^{ss}(p)$ is $X^t$-invariant only for $t \in
\Pi_X(p) \cdot \Z$.
The \emph{index} of $\gamma$ is the dimension of the stable
manifold $W^s(\gamma) = W^s(p)$.

\begin{lemma} \label{l.irrational}
 Given any $r \in \NN$, there exists a residual subset $\RR_1$ of
 $\xisr$ such that if $X \in \RR_1$ then given
 any distinct closed orbits $\gamma$, $\gamma'$, we have that
 $$
 \frac{\Pi_X(\gamma)}{\Pi_X(\gamma')} \in \R \setminus \Q.
 $$
 \end{lemma}

 \begin{proof}
 For $N \in \N$, let $A_N \subset \xisr$ be
 the set
 of vector fields $X$ such that all singularities of
 $X$
 are hyperbolic and all closed orbits with periods
 less than $N$ are
 hyperbolic. It follows from the standard proof of
 the Kupka--Smale
 theorem that the set $A_N$ is open and dense in
 $\xisr$.

 Now let $a_1$, $a_2$, \dots be an enumeration of the
 positive rational numbers
 and let $B_N \subset \xisr$ be the set of vector
 fields $X \in A_N$
 such that if $\gamma$, $\gamma'$ are distinct closed
 orbits with periods
 less than $N$ then $\Pi_X(\gamma) / \Pi_X(\gamma')$
 does not belong to
 $\{a_1, \dots, a_N \}$.

 If $X \in A_N$ then the number of orbits with
 periods less that $N$ is finite.
 Moreover, each of these orbits has a continuation
 and the period varies continuously.
 It follows that the set $B_N$ is open.

 Let us show that $B_N$ is also dense, so
 we can define $\RR_1 = \cap_N B_N$.
 Given any $X_0 \in \xis$, first approximate it by
 $X_1 \in A_N$.
 Let $\gamma_1$, \dots, $\gamma_k$ be the
 $X_1$-orbits with periods less than $N$.
 Let $\delta>0$ be small enough such that
 the neighborhoods $B(\gamma_i, \delta)$ are disjoint.
 Take $C^r$ functions $\psi_i : M \to [0,1]$ such
 that
 $\psi_i$ equals $1$ in $\gamma_i$ and equals $0$
 outside $B(\gamma_i, \delta)$.
 For $s \in \R_{+}^k$ close to $0$, let
 $$
 Y_s = \left( \prod_{i=1}^k (1+s_i \psi_i)^{-1}
 \right) X_1
 $$
 Then $Y_s$ has the same orbits as $X_1$ and $Y_s$
 converges to $X_1$
 in the $C^r$ topology as $s \to 0$.
 Moreover,
 $$
 \Pi_{Y_s}(\gamma_i) = (1+s_i) \Pi_{X_1}(\gamma_i),
 $$
 so we can find $s\in \R_{+}^k$ close to $0$ such
 that
 $Y_s \in A_N$ and
 the quotients $\Pi_{Y_s}(\gamma_i) /
 \Pi_{Y_s}(\gamma_j)$, $i\neq j$,
 do not intersect $\{a_1, \ldots, a_N\}$.
 If $\gamma$ is another closed orbit of $Y_s$, then
 $\Pi_{Y_s}(\gamma) \geq \Pi_{X_1}(\gamma) \geq N$.
 This proves that $Y_s \in B_N$.
 \end{proof}

We shall use the following simple fact, whose proof is omitted (it
follows easily from the transitivity of the future orbits of
irrational rotations of the circle):

\begin{lemma} \label{l.elementary}
Given numbers $a>0$, $b>0$ and $\eps>0$, with $a/b$ irrational,
the set $$ \big\{ ma + nb + s: \; m, n \in \N, \ |s| < \eps \big\}
$$ contains an interval of the form $[T, +\infty)$.
\end{lemma}

\medskip

We may now prove Theorem B:

\begin{proof}[Proof of Theorem B]
Let $\RR_1 \subset \xisr$ be the residual set given by
Lemma~\ref{l.irrational}, and let $H = H_X(\gamma_0)$ be a
non-trivial homoclinic class of some $f \in \RR_1$. Take two
nonempty open sets $U$, $V$ intersecting $H$. We shall prove that
there exists $t_0 > 0$ such that $X^t(U) \cap V \neq \varnothing$
for every $t \geq t_0$.

Let $\gamma$ and $\gamma'$ be distinct periodic orbits in $H$ with
same index, such that $\gamma  \cap U \neq \varnothing$ and
$\gamma'\cap V \neq \varnothing$.
In order to simplify the notation, let $a = \Pi_X(\gamma)$ and $b
= \Pi_X(\gamma')$. Recall that $a/b \in \R \setminus \Q$.

Take $p \in \gamma \cap U$ and $q\in \gamma' \cap V$. Notice that
$W^u (p) \cap W^s(q)$ is non-empty. Fix a point $y$ in this
intersection. There exists $\tau_1 > 0$ such that $X^{-\tau_1}(y)
\in W^{uu}(p)$ and, consequently, the sequence $\{X^{(- \tau_1 -
ma)}(y)\}_{m\in \N}$ is contained in $W^{uu}(p)$ and converges to
$p$. Therefore we can find $t_1 > 0$ such that $$ X^{(- t_1 - ma)}
(y) \in U  \text{ for every $m \in \N$.} $$ Analogously, there
exist $t_2 >0$ and $\eps > 0 $ such that $$ X^{(  t_2 + nb +
s)}(y) \in V \text{ for every $m \in \N$ and $|s| < \eps$.} $$ Let
$T>0$, depending on $a$, $b$ and $\eps$, be given by
Lemma~\ref{l.elementary}. Set $t_0 = t_1 + t_2 + T$. Then, for any
$t \geq t_0$, there exist numbers $m$, $n\in\N$ and $|s|<\eps$
such that $t = t_1 + t_2 + ma + nb + s$. So $X^t(U) \cap V$
contains the point $X^{(t_2 + nb + s)}(y)$. This concludes the
proof.
\end{proof}

Now we explain how Theorem A follows from Theorem B. We first need
a couple of definitions:

\begin{definition}
Let $\Lambda$ be a compact invariant set of $X \in \xis$. Then we
set $P_i(\Lambda) \equiv \{ p \in \Lambda: p \text{ is a
hyperbolic periodic point of $X$ with index $i$} \}$.
\end{definition}

The next definition comes from~\cite{BDP}:

\begin{definition}
Let $p$ be a periodic point of a flow $X \in \xis$ and $U$ be a
neighborhood of $p$ in $M$. Then the \emph{homoclinic class of $p$
relative to $U$} is given by $$ \HR_X(p, U) \equiv \mathrm{cl} \,
\{q \in H_X(p) \cap \text{Per}(X): \text{the orbit } \gamma(q)
\text{ is contained in } U \}. $$ It is easily seen that $\HR_X(p,
U)$ is a compact transitive invariant set. Moreover, if
$\text{ind}(p) = i$, then $P_i(\HR_X(p, U))$ is dense in $\HR_X(p,
U)$.
\end{definition}

We need the following lemma, which is a consequence of a theorem
by Arnaud~\cite{Ar} together with an argument from~\cite{BDP}:

\begin{lemma} \label{arnaud}
There is a residual subset $\RR$ of $\xis$ such that if $\Lambda$
is an isolated transitive set of $X \in \RR$, then $\Lambda =
\HR_X(p, U)$ for some periodic point $p \in \Lambda$.
\end{lemma}

\begin{proof}
Let $\RR_2$ be as in Theorem 1 of~\cite{Ar} and let $\RR_3$ be as
in Theorem B of~\cite{BD2}, and set $\RR \equiv \RR_2 \cap \RR_3$.
Let $\Lambda$ be an isolated transitive set of $X \in \RR$, with
$U$ an isolating block of $\Lambda$ .

By Theorem 1 of~\cite{Ar}, there is a sequence of periodic orbits
$\gamma_k$ which converge to $\Lambda$ in the Hausdorff topology.
The orbit $\gamma_k$ is contained in $U$ for $k$ sufficiently
large. Since $\Lambda$ is the maximal invariant set of $U$, it
follows that for large $k$ the orbit $\gamma_k$ is contained in
$\Lambda$. Since the sequence $\{\gamma_k\}$ converges to
$\Lambda$ in the Hausdorff topology, the set of periodic points
contained in $\Lambda$ must be a dense subset of $\Lambda$.

Now, since $\Lambda$ is transitive and has a dense subset of
periodic points, we apply an argument from [BDP] which uses
Theorem B of~\cite{BD2} to conclude that given any periodic point
$p \in \Lambda$ then $$ \Lambda = \HR_X(p, U).$$

\end{proof}

We are now ready to prove Theorem A:

\begin{proof}[Proof of Theorem A]
It is easy to see that the proof of Theorem B actually implies the
following result:

\begin{theorem D}
There is a residual subset $\RR$ of $\xis$ such that if $Y \in
\RR$ and $\Lambda$ is a non-trivial transitive set of $Y$ such
that for some $i \in \{1, \ldots, d-1\}$ the set $P_i(\Lambda)$ is
dense in $\Lambda$, then $\Lambda$ is topologically mixing for
$Y$.
\end{theorem D}

Now by Lemma \ref{arnaud} above we have that $\Lambda$ coincides
with some relative homoclinic class $\HR_X(p, U)$. Let $i$ be the
index of the periodic point $p$. Since $P_i(\HR_X(p, U))$ is dense
in $\HR_X(p, U)$, we conclude that $\Lambda$ satisfies the
hypotheses of Theorem D above, and therefore that $\Lambda$ is a
mixing set for $X$.
\end{proof}

At last, we give the:

\begin{proof}[Proof of Theorem C]
Let $S$ be a compact $3$-manifold and let $\diff$ be the set of
$C^1$ diffeomorphisms of $S$ endowed with the $C^1$-topology. The
key of the construction is the following result of Bonatti and
D{\'\i}az (\cite[Theorem 3.2]{BD3}): \emph{ There exist an open
set $\mathcal{U}_0 \subset \diff$ and a dense subset
$\mathcal{D}_0 \subset \mathcal{U}_0$ such that for every $f \in
\mathcal{D}_0$ there are an open set $B \subset S$ and an integer
$n \in \N$ such that every $x \in B$ is a periodic point of $f$ of
(prime) period $n$. }

Let $f_0: S \to S$ be a diffeomorphism from the set
$\mathcal{D}_0$ above. Let $X_0^t: M \to M$ be the suspension
flow. As usual, $M$ is the $4$-manifold obtained from $S \times
[0,1]$ by gluing points $(x,1)$ and $(f_0(x), 0)$. We will
identify $S$ with the submanifold $\{(x,0)\in M; \; x \in S\}$ of
$M$.

Let $\mathcal{U} \subset \xis$ be a small neighborhood of $X_0$
such that every  vector field $X \in \mathcal{U}$ is transverse to
$S$ and, moreover, the first-return map $f_X \in \diff$ belongs to
$\mathcal{U}_0$.  For $X \in \mathcal {U}$, we let $\tau_X:S \to
\R_+$ be the return-time map, which is a $C^1$-smooth function
depending continuously (in the $C^1$ topology) on $X \in \mathcal
{U}$.

We will omit the proof of the following:
\begin{lemma} \label{l.perturbed return}
For every $X \in \mathcal{U}$ and every neighborhood $\mathcal{V}
\ni X$, if $\tilde{f}$ is a small perturbation of $f_X$ and
$\tilde \tau$ is a small perturbation of $\tau_X$ then there is
$\tilde{X} \in \mathcal{V}$ such that $f_{\tilde{X}} = \tilde{f}$
and $\tau_{\tilde X}=\tilde \tau$.
\end{lemma}

Now let $X_1 \in \mathcal{U}$. We shall prove that there exists
$X_4$ arbitrarily close to $X_1$ which has a non-trivial
homoclinic class which is not topologically mixing.

Let $f_1 = f_{X_1}$, and $\tau_1=\tau_{X_1}$. Take $f_2 \in
\mathcal{D}_0$ close to $f_1$. Since $f_2 \in \mathcal{D}_0$,
there is a ball $B \subset S$ of points that are $f_2$-periodic,
of period $n$. Let $\tau_1^n:S \to \R_+$ be defined by
$\tau_1^n=\sum_{j=0}^{n-1} \tau_1 \circ f_1^j$.

Using a chart, we identify $B$ with a ball $B(0,r)\subset \R^3$,
in such a way that the kernel of the differential $D\tau_1^n(0)$
contains the plane $xy$.

Let $f_3\in \diff$ be a perturbation of $f_2$ such that:
\begin{itemize}
\item $f_3$ equals $f_2$ outside $\cup_{j=0}^{n-1} f_2^j(B)$;
\item there exists a ball $B_1 = B(0,r_1)$, with $0<r_1<r$, such that
$f_3^n(B_1) = B_1$;
\item $f_3^n$ restricted to $B_1$ is a orthogonal
rotation (indicated by $R$) of angle $2 \pi /m$, where $m\in \N$,
along the axis $y$.
\end{itemize}

It is easy to construct a map $\tau_3:S \to \R$ $C^1$ close to
$\tau_1$, such that $\tau_3^n=\sum_{j=0}^{n-1} \tau_3 \circ f_3^j$
is an affine map in a smaller ball $B_2$ around $0$ and such that
$D\tau_3^n(0)=D\tau_1^n(0)$. That is, if $x \in B_2$ then
$\tau_3^n(x) = \tau_3^n(0) + D\tau_3^n(0) \cdot x$.

Using Lemma~\ref{l.perturbed return}, we find a flow $X_3$ close
to $X_1$ and such that $f_{X_3}=f_3$ and $\tau_{X_3}=\tau_3$.

Let $x\in B_2 \setminus\{0\}$. Its successive returns to $B_1$
under the flow $X_3$ are $R(x)$,\ldots,$R^{m-1}(x)$, $R^m(x) = x$.
In particular, $x$ is a periodic point. Summing the respective
return times we get that the period of $x$ is $\sum_{j=0}^{m-1}
\tau_3^n(R^j(x)) =m \tau_3^n(0)$, since $\sum_{j=0}^{m-1} R^j(x)$
belongs to the $y$ axis. That is, all points in $B_2
\setminus\{0\}$ are periodic under $X_3$ of (prime) period $m
\tau_3^n(0)$.

Let $B_3 \subset B_2$ be a ball (not centered in $0$) such that
cl\,$B_3$, cl\,$R(B_3)$,\ldots, cl\,$R^{m-1}(B_3)$ are pairwise
disjoint.

Choose now some $f_4:S \to S$ which is $C^1$ close to $f_3$ and
such that:
\begin{itemize}
\item $f_4$ equals $f_3$ outside $B_3$;
\item $f_4^{nm}$ restricted to $B_3$ has a non-trivial homoclinic
class (say, a solenoid attractor).
\end{itemize}

Let also $\tau_4:S \to \R_+$ be given by $\tau_4=\tau_3 \circ
f_3^{-1} \circ f_4$.  Then $\tau_4$ is $C^1$ close to $\tau_3$.
Using Lemma~\ref {l.perturbed return} again, we obtain $X_4$ $C^1$
close to $X_3$ such that $f_{X_4} = f_4$ and $\tau_{X_4}=\tau_4$.
The return time of a point $x \in B_3$ to $B_3$ under $X_4$ is
$\tau_4^{nm}(x) = \tau_3^{nm}(f_3^{-1} \circ f_4 (x))$, which
independs of $x$ (where, as usual, we let
$\tau_4^{nm}=\sum_{j=0}^{nm-1} \tau_4 \circ f_4^j$ and
$\tau_3^{nm}=\sum_{j=0}^{nm-1} \tau_3 \circ f_3^j$). Therefore
$X_4^t$ has a non-trivial homoclinic class which is not
topologically mixing.
\end{proof}

{\bf Acknowledgements:}  We would like to thank the referee for
his detailed report and for pointing out to us that our proof of
Theorem~B was not restricted to the $C^1$ topology.


\bigskip

\small
\begin{flushleft}
Flavio Abdenur ({\tt flavio@impa.br})

IMPA

Estr. D. Castorina 110

22460-320 Rio de Janeiro -- Brazil.

\bigskip

Artur Avila ({\tt avila@impa.br})

Coll\`ege de France

3 rue d'Ulm

75005 Paris -- France.

\bigskip

Jairo Bochi ({\tt bochi@impa.br})

IMPA

Estr. D. Castorina 110

22460-320 Rio de Janeiro -- Brazil.

\end{flushleft}

\end{document}